\documentclass[11pt,openbib]{article} 
\usepackage{amsfonts,amsthm}
\usepackage{amssymb,amsmath}
\newcommand{\R}{\mathbb{R}} 
\newcommand{\C}{\mathbb{C}}
 
\newcommand{\Z}{\mathbb{Z}}
 
\newcommand{\gl}{\mathop{\mathrm{gl}}\nolimits}

\newcommand{\ad}{\mathop{\mathrm{ad}}\nolimits}

\newcommand{\comp}{\raisebox{0pt}{$\scriptstyle\circ \, $}}

\newcommand{\setrule}{\, \mathop{\rule[-4pt]{.5pt}{13pt}\, }\nolimits}
\newcommand{\smallspace}{\smallskip\par\noindent}

\newcommand{\spann}{\mathop{\rm span}\nolimits} 

\allowdisplaybreaks
\begin{document}
\begin{center}
{\large \bf The Weyl group of a fundamental sandwich algebra}\footnotetext{version: \today} \\
\vspace{.05in} 
Richard Cushman\footnotemark 
\end{center}
\footnotetext{\parbox[t]{3in}{Department of Mathematics and Statistics, \\ 
University of Calgary, 
Calgary, AB, T2N 1N4 Canada} }

In this paper we construct a Weyl group $W_{\widetilde{\mathcal{R}}}$ for a fundamental sandwich 
algebra $\widetilde{\mathfrak{g}}$, whose definition and basic properties we now recall. \medskip  

$\widetilde{\mathfrak{g}} = \mathfrak{g} \oplus \widetilde{\mathfrak{n}}$ is a very special sandwich algebra 
provided that 1) $\widetilde{\mathfrak{g}}$ is a subalgebra of the complex simple Lie algebra 
$\underline{\mathfrak{g}}$ of rank $1$ greater than the rank of the simple Lie algebra $\mathfrak{g}$ with Cartan subalgebra $\mathfrak{h}$; 2) $\widetilde{\mathfrak{n}}$ is a sandwich, that is, 
$[\widetilde{\mathfrak{n}}, [\widetilde{\mathfrak{n}}, \widetilde{\mathfrak{n}}]] =0$ and 
$[\widetilde{\mathfrak{n}}, \widetilde{\mathfrak{n}}] \ne 0$, which is 
the nilpotent radical of $\widetilde{\mathfrak{g}}$; 3) ${\ad }_{\mathfrak{h}}$ is a maximal torus of 
$\gl (\widetilde{\mathfrak{n}}, \C )$. A very special 
sandwich algebra $\widetilde{\mathfrak{g}}$ is fundamental if and only if the center $Z$ of the nilradical 
$\widetilde{\mathfrak{n}}$ is $1$-dimensional. Since ${\ad }_{\mathfrak{h}}$ is a maximal torus, we may write 
$\widetilde{\mathfrak{n}} = \sum_{\widehat{\alpha } \in \mathfrak{R} } \oplus 
{\widehat{\mathfrak{g}}}_{\widehat{\alpha }}$, where 
${\widehat{\mathfrak{g}}}_{\widehat{\alpha }}$ is a $1$-dimensional root space spanned by the nonzero root vector $X_{\widehat{\alpha }}$ for every root $\widehat{\alpha } \in \mathfrak{R}$. Let 
$\zeta $ be the root in $\mathfrak{R}$ such that $Z = {\spann }_{\C } \{ X_{\zeta } \} $. Because 
$\widetilde{\mathfrak{g}}$ is fundamental, it follows that 
$\zeta |\mathfrak{h} = \widehat{0}$ is the linear function on $\mathfrak{h}$ which is identically zero. We have 
$\widehat{\mathcal{R}} = \mathfrak{R} \setminus \{ \zeta  \} = \widehat{\Pi }\, \amalg (-\widehat{\Pi })$, where 
$\widehat{\Pi } = {\{ {\widehat{\alpha }}_i \} }^M_{i=1}$ is a set of positive roots. In \cite{cushman17b} we have shown that $\widehat{\mathcal{R}}$ is a system of roots, which we define below. 

\section{Weyl group of a system of roots $\widehat{\mathcal{R}}$}

We begin by constructing a Weyl group associated to a system of roots 
$\widehat{\mathcal{R}}$ of the fundamental sandwich algebra $\widetilde{\mathfrak{g}}$. \medskip 

We recall the definition of a system of roots. Let $V$ be a finite dimensional real vector space with $\Phi$ a finite subset of nonzero vectors. 
Recall that a system of roots $\Phi $ satisfies the axioms: \medskip 

\noindent 1. Let $V = {\spann }_{\R }\Phi $, using addition $+$ of vectors in $V$.
\smallspace 
\noindent 2. $\Phi = - \Phi $, where $-$ is the additive inverse of $+$. 
\smallspace
\noindent 3. \parbox[t]{4.75in}{For every $\beta $, $\alpha \in \Phi $ there is an \emph{extremal root chain 
${\mathcal{S}}^{\beta }_{\alpha }$ through $\beta $ in the direction $\alpha $} given by 
$\{ \beta +j \alpha \in \Phi \cup \{ 0 \} \setrule \mbox{for every $j \in \Z$, $-q \le j \le p $} \} $. Here 
$q,p \in {\Z }_{\ge 0}$ and are as large as possible. The pair $(q,p)$ is the \linebreak 
\emph{integer pair} associated to ${\mathcal{S}}^{\beta }_{\alpha }$. The integer $\langle \beta , \alpha \rangle = q-p$ 
is called the \emph{Killing integer} of ${\mathcal{S}}^{\beta }_{\alpha }$.} 
\smallspace
\noindent  4. Fix $\alpha \in \Phi $ and suppose that ${\beta }_1$, ${\beta }_2$, and 
${\beta }_1+ {\beta }_2 \in \Phi $. Then 
\begin{equation}
\langle {\beta }_1+{\beta }_2, \alpha \rangle = \langle {\beta }_1, \alpha \rangle + \langle {\beta }_2, \alpha \rangle .
\label{eq-sec1subsec1one}
\end{equation}
\noindent 5. For every $\alpha \in \Phi $ we have $\langle \alpha , \alpha \rangle = 2$.  

\subsection{Definition of the Weyl group $W_{\widehat{\mathcal{R}}}$}

Consider the system of roots $\widehat{\mathcal{R}}$ associated to the fundamental sandwich algebra 
$\widetilde{\mathfrak{g}}$. \medskip 

We need the notion of a reflection before we can define the Weyl group $W_{\widehat{\mathcal{R}}}$ of the system of roots $\widehat{\mathcal{R}}$. From axiom $4$ of a system of roots it follows that for every fixed $\alpha \in \widehat{\mathcal{R}} $, the function $K_{\alpha }: \widehat{\mathcal{R}} \subseteq V \rightarrow \Z: \beta \mapsto \langle \beta , \alpha \rangle $ 
is linear, that is, if $\gamma $, $\delta \in \widehat{\mathcal{R}} $ and $\gamma + \delta \in \widehat{\mathcal{R}}$, then 
$K_{\alpha }(\gamma +\delta ) = K_{\alpha }(\gamma ) + K_{\alpha }(\delta )$. Knowing the values of 
$K_{\alpha }$ on a basis $B \subseteq \widehat{\mathcal{R}} $ of the real vector space $V$, we can extend $K_{\alpha }$ uniquely to a real linear function $k_{\alpha }: V \rightarrow \R $ such that $k_{\alpha }|B = K_{\alpha }|B$. Because $k_{\alpha }$ is unique, we have $k_{\alpha }|\widehat{\mathcal{R}} = K_{\alpha }$. 
Since $k_{\alpha }(\alpha ) = 2$ by axiom $5$, it follows that $k_{\alpha }$ is nonzero. 
For each root $\alpha \in \widehat{\mathcal{R}} $ the real linear mapping 
\begin{equation}
{\sigma }_{\alpha }: V \rightarrow V: v \mapsto v -k_{\alpha }(v)\alpha 
\label{eq-onestar}
\end{equation}
is a \emph{reflection} in the hyperplane $H_{\alpha } = \{ v \in V \setrule k_{\alpha }(v) = 0 \} $, since 
${\sigma }_{\alpha }$ leaves every vector in $H_{\alpha }$ fixed and ${\sigma }_{\alpha }(\alpha ) = 
-\alpha $, because $k_{\alpha }(\alpha ) =2$. Moreover, the line ${\spann }_{\R } \{ \alpha \} $ is 
transverse to the hyperplane $H_{\alpha }$.  \medskip  

\noindent \textbf{Lemma 1.1.1} For every $\alpha \in \widehat{\mathcal{R} }$ the real linear mapping 
${\sigma }_{\alpha }$ (\ref{eq-onestar}) is an involution, that is, ${\sigma }_{\alpha } \comp {\sigma }_{\alpha} 
= {\mathrm{id}}_V$. \medskip 

\noindent \textbf{Proof.} For every $v \in V$ we have 
\begin{align}
{\sigma }_{\alpha }\big( {\sigma }_{\alpha }(v) \big) & = {\sigma }_{\alpha }(v) - k_{\alpha }({\sigma }_{\alpha }(v))\alpha 
\notag \\
& = v - k_{\alpha }(v)\alpha  - k_{\alpha } \big( v -k_{\alpha }(v) \alpha \big) \alpha \notag \\
& = v - k_{\alpha }(v)\alpha  -k_{\alpha }(v)\alpha  + k_{\alpha } (\alpha ) k_{\alpha }(v) \alpha , \notag \\
& \hspace{.5in}\mbox{since $k_{\alpha }$ is a real linear function on $V$} \notag \\
& = v, \, \, \, \mbox{since $k_{\alpha }(\alpha ) =2$ by axiom $5$.} \tag*{$\square $}
\end{align} 

\noindent \textbf{Corollary 1.1.1a} For every $\alpha \in \widehat{\mathcal{R}}$, the reflection 
${\sigma }_{\alpha }$ (\ref{eq-onestar}) is an invertible 
real linear mapping of $V$ into itself. \medskip 

\noindent \textbf{Corollary 1.1.1b} For every $\alpha \in \widehat{\mathcal{R}}$ the reflection 
${\sigma }_{\alpha }$ (\ref{eq-onestar}) sends $\widehat{\mathcal{R}} \subseteq V$ into itself. \medskip 

\noindent \textbf{Proof}. Suppose that $\beta \in \widehat{\mathcal{R}}$. Then 
${\sigma }_{\alpha }(\beta ) = \beta -\langle \beta , \alpha \rangle \alpha = \beta +j \alpha $, 
where $j = -\langle \beta , \alpha \rangle = p -q$. Look at the extremal root chain 
${\mathcal{S}}^{\beta }_{\alpha}$
\begin{displaymath}
\beta - q \alpha , \ldots , \beta - \alpha , \beta , \beta + \alpha , \ldots , \beta + p \alpha , 
\end{displaymath}
which has length $p+q+1$. Counting back $q$ nodes on ${\mathcal{S}}^{\beta }_{\alpha}$ from 
the node $\beta + p \alpha $ gives $\beta + (p-q) \alpha $, which is a node on ${\mathcal{S}}^{\beta }_{\alpha}$ 
since $p+q+1 >q$. Thus ${\sigma }_{\alpha }(\beta ) \in \widehat{\mathcal{R}} \cup \{ \zeta \}$. 
Since $\beta \ne \zeta $ and ${\sigma }_{\alpha }$ is invertible by corollary 1.1.1a, it follows that 
${\sigma }_{\alpha }(\beta ) \ne \zeta $. Hence ${\sigma }_{\alpha }(\beta ) \in \widehat{\mathcal{R}}$. 
\hfill $\square $ \medskip

\noindent \textbf{Lemma 1.1.2} For every $\alpha $, $\beta \in \widehat{\mathcal{R}} $ the reflection 
${\sigma }_{\alpha }$ maps the extremal 
root chain ${\mathcal{S}}^{\beta }_{\alpha }$ with integer pair $(q,p)$ into itself. In particular, for 
every $j \in \Z $ with $-q \le j \le p$ we have ${\sigma }_{\alpha }(\beta + j \alpha ) = 
\beta + \big( p-(q+j) \big) \alpha $. \medskip 

\noindent \textbf{Proof.} We compute
\begin{align}
{\sigma }_{\alpha }(\beta + j \alpha ) & = \beta +j\alpha  - k_{\alpha }(\beta +j\alpha ) \alpha \notag \\
&\hspace{-.5in}  = \beta +j \alpha -k_{\alpha}(\beta )\alpha -j k_{\alpha }(\alpha ) \alpha , \, \, \, \mbox{since $k_{\alpha }$ is linear} \notag \\
& \hspace{-.5in} = \beta + \big( p - (q+j) \big) \alpha , \, \, \, \mbox{since $k_{\alpha }(\alpha ) =2$ and $-k_{\alpha}(\beta ) = p-q$.}
\tag*{$\square $} 
\end{align} 

\noindent \textbf{Corollary 1.1.2a} For every $\alpha \in \widehat{\mathcal{R}} $ the map  
${\sigma }_{\alpha }$ is bijective on $\widehat{\mathcal{R}} $.  \medskip 

\noindent \textbf{Proof.} From the lemma it follows that ${\sigma }_{\alpha }$ maps $\widehat{\mathcal{R}}$ into itself. Because ${\sigma }_{\alpha }$ is an invertible linear mapping to $V$ onto itself, ${\sigma }_{\alpha }$ maps only $0$ onto $0$. Hence ${\sigma }_{\alpha }$ maps $\widehat{\mathcal{R}} $ bijectively to itself. \hfill $\square $ \medskip 

Since the set of roots $\widehat{\mathcal{R}}$ is finite by hypothesis, for every $\alpha \in \widehat{\mathcal{R}} $ the reflection 
${\sigma }_{\alpha }$ is a permutation of the elements of $\widehat{\mathcal{R}} $. Let $W_{\widehat{\mathcal{R}}}$ be the group 
generated by the reflections ${\sigma }_{\alpha }$ with $\alpha \in \widehat{\mathcal{R}} $. We call 
$W_{\widehat{\mathcal{R}}}$ the \emph{Weyl group} associated to the system of roots 
$\widehat{\mathcal{R}}$. Since $W_{\widehat{\mathcal{R}}}$ is 
a subgroup of the group of permutations of the elements of the finite set $\widehat{\mathcal{R}} $, it is a finite group. 

\subsection{Generators and relations for $W_{\widehat{\mathcal{R}}}$}

In theorem 3.4.6 of \cite{cushman17b} we have shown that $\widehat{\mathcal{R}} = \mathfrak{R} \setminus 
\{ \zeta \}$ is a system of roots, where $\mathfrak{R}$ is the set of roots of the 
nilradical $\widetilde{\mathfrak{n}}$ of a fundamental sandwich algebra. 
From now on we deal only with such a system of roots. \medskip 

Next we determine the generators and relations of the Weyl group $W_{\widehat{\mathcal{R}}}$. 
By definition the reflections ${\sigma }_{\widehat{\alpha }}$, $\widehat{\alpha } \in \widehat{\mathcal{R}} $ generate the Weyl group $W_{\widehat{\mathcal{R}}}$. Claim 1.2.2 will sharpen this. First we show \medskip 

\noindent \textbf{Lemma 1.2.1} For every $\widehat{\alpha } \in \widehat{\mathcal{R}} $ we have 
${\sigma }_{\widehat{\alpha }} = ({\sigma }_{-\widehat{\alpha }})^{-1}$. \medskip 

\noindent \textbf{Proof.} Let $v \in V = {\spann }_{\R } \widehat{\mathcal{R}} $. For every 
$\widehat{\alpha } \in \widehat{\mathcal{R}} $ we compute
\begin{align}
({\sigma }_{-\widehat{\alpha }} \comp {\sigma }_{\widehat{\alpha }})(v) & = 
{\sigma }_{\widehat{\alpha }}(v) - k_{-\widehat{\alpha }}\big( {\sigma }_{\widehat{\alpha }}(v) \big) (-\widehat{\alpha }),   
\, \, \, \mbox{since $-\widehat{\alpha } \in \widehat{\mathcal{R}}$} \notag \\
& = v - k_{\widehat{\alpha }}(v)\widehat{\alpha } +k_{-\widehat{\alpha }}(v-k_{\widehat{\alpha }}(v)\widehat{\alpha })\widehat{\alpha } 
\notag \\
& = v - k_{\widehat{\alpha }}(v)\widehat{\alpha } +k_{-\widehat{\alpha }}(v) \widehat{\alpha } 
-k_{\widehat{\alpha }}(v)k_{-\widehat{\alpha }}(\widehat{\alpha}) \widehat{\alpha } .  \notag 
\end{align}  
Using lemma 1.2.3 of \cite{cushman17b} we get $k_{-\widehat{\alpha }}(\beta ) = \langle \beta, -\widehat{\alpha } \rangle = - \langle \beta , \widehat{\alpha } \rangle = -k_{\widehat{\alpha }}(\beta)$ 
for every $\beta \in \widehat{\mathcal{R}}$. Consequently, 
$k_{-\widehat{\alpha }}(v) = - k_{\widehat{\alpha }}(v)$ for every $v \in V$, since 
$B \subseteq \widehat{\mathcal{R}} $ is a basis of $V$. So 
\begin{align}
({\sigma }_{-\widehat{\alpha }} \comp {\sigma }_{\widehat{\alpha }})(v) & = v - k_{\widehat{\alpha }}(v) \widehat{\alpha } 
-k_{\widehat{\alpha }}(v)\widehat{\alpha } + k_{\widehat{\alpha }}(\widehat{\alpha })k_{\widehat{\alpha}}(v) \widehat{\alpha } 
\notag \\
& = v, \, \, \, \mbox{since $k_{\widehat{\alpha }}(\widehat{\alpha }) =2$.} \tag*{$\square $} 
\end{align} 

Let $\widehat{\Pi } = {\{ {\widehat{\alpha }}_i \} }^M_{i=1}$ be the set of positive roots in 
$\widehat{\mathcal{R}} = { \{ \pm \widehat{\alpha } \} }^M_{i=1} = { \{ {\widehat{\beta }}_i \} }^{2M}_{i=1}$, see 
claim 2.1.6 of \cite{cushman17b}. We now prove \medskip 

\noindent \textbf{Claim 1.2.2} $W_{\widehat{\mathcal{R}}}$ is generated by the 
reflections ${\sigma }_{\widehat{\alpha }}$ with $\widehat{\alpha} \in \widehat{\Pi }$.  \medskip 

\noindent \textbf{Proof.} Let $w \in W_{\widehat{\mathcal{R}}}$. Then $w = {\sigma }_{{\widehat{\beta }}_{i_1}} \comp 
\cdots {\sigma }_{{\widehat{\beta }}_{i_n}} $ where ${\widehat{\beta }}_{i_j} \in \widehat{\mathcal{R}} $ for all 
$1 \le j \le n$ and 
$i_j \in \{ 1, \dots , 2M \}$. Now $\widehat{\mathcal{R}} = \widehat{\Pi } \, \amalg (-\widehat{\Pi })$. For every 
$i_{\ell } \in \{ i_1, \ldots , i_n \} $ such that ${\widehat{\beta }}_{i_{\ell }} \in (-\widehat{\Pi })$, we have 
${\sigma }_{{\widehat{\beta }}_{i_{\ell }}} = ({\sigma }_{-{\widehat{\beta }}_{i_{\ell }}})^{-1} = 
{\sigma }_{-{\widehat{\beta }}_{i_{\ell }}}$, because ${\sigma }_{{\widehat{\beta }}_{i_{\ell }}}$ is a reflection. 
Note that $-{\widehat{\beta }}_{i_{\ell }} \in \widehat{\Pi }$. Thus $w$ is a composition of reflections 
${\sigma }_{\widehat{\alpha }}$, where $\widehat{\alpha } \in \widehat{\Pi }$. Hence 
$W_{\widehat{\mathcal{R}}}$ is generated 
by the reflections ${\sigma }_{\widehat{\alpha }}$, where $\widehat{\alpha } \in \widehat{\Pi }$. \hfill $\square $ \medskip 

\noindent \textbf{Claim 1.2.3} The Weyl group $W_{\widehat{\mathcal{R}}}$ is an abelian group 
generated by the reflections ${\sigma }_{{\widehat{\alpha }}_i}$, $i \in J = \{ 1, \ldots , M \} $, which satisfy the relations 
\begin{displaymath}
\begin{array}{lrl}
1. & {\sigma }^2_{{\widehat{\alpha }}_i} & = {\sigma }_e, \, \,\mbox{for all $i \in J$} \notag \\
2. & {\sigma }_{{\widehat{\alpha }}_i} \comp {\sigma }_{{\widehat{\alpha }}_j} & = 
{\sigma }_{{\widehat{\alpha }}_j} \comp {\sigma }_{{\widehat{\alpha }}_i}, \, \, \mbox{for all $i,j \in J$} \notag \\
3. & {\sigma }_{{\widehat{\alpha }}_1} \comp {\sigma }_{{\widehat{\alpha }}_2} \cdots \comp  
{\sigma }_{{\widehat{\alpha }}_M} & = - {\sigma }_e. 
\end{array}
\end{displaymath}
Here ${\sigma }_e = {\mathrm{id}}_{V}$. \medskip 

\noindent \textbf{Proof.} The claim is a consequence of the following lemmas. \medskip 
 
\noindent \textbf{Lemma 1.2.3a} For every $i,j \in J = \{ 1, 2, \ldots , M \} $ with 
$i \ne j$ we have $\langle {\widehat{\alpha }}_j, {\widehat{\alpha }}_i \rangle =0$. \medskip 

\noindent \textbf{Proof.} Suppose that ${\widehat{\alpha }}_j + 
{\widehat{\alpha }}_i \in \widehat{\Pi }$, then 
$0 \ne X_{{\widehat{\alpha}}_1+{\widehat{\alpha }}_2} = [X_{{\widehat{\alpha }}_j}, X_{{\widehat{\alpha }}_i}] \in Z = {\widehat{\mathfrak{g}}}_{\widehat{0}}$. So for some nonzero complex number $c$ we have $[X_{{\widehat{\alpha }}_j}, X_{{\widehat{\alpha }}_i}] 
= c\, X_{\widehat{0}}$. This  implies that ${\widehat{\alpha }}_j + {\widehat{\alpha }}_i = \widehat{0}$. Thus  
${\widehat{\alpha }}_j = - {\widehat{\alpha }}_i \notin \widehat{\Pi }$, which 
contradicts the hypothesis that ${\widehat{\alpha }}_j \in \widehat{\Pi }$. Therefore 
${\widehat{\alpha }}_j + {\widehat{\alpha }}_i \notin \widehat{\Pi }$ and hence does not lie in 
$\widehat{\mathcal{R}} \cup \{ \widehat{0} \} $. Now suppose that ${\widehat{\alpha }}_j - 
{\widehat{\alpha }}_i \in \widehat{\mathcal{R}} \cup \{ \widehat{0} \} $. Then an argument similar to the one given above shows that ${\widehat{\alpha }}_j -{\widehat{\alpha }}_i = \widehat{0}$, that is, ${\widehat{\alpha }}_j = {\widehat{\alpha }}_i$. But this contradicts the hypothesis that $j \ne i$ and thus 
${\widehat{\alpha }}_j \ne {\widehat{\alpha }}_i$. Therefore the root chain 
${\mathcal{S}}^{{\widehat{\alpha }}_j}_{{\widehat{\alpha }}_i}$ in the system of roots 
$\widehat{\mathcal{R}}$ is extremal with integer pair $(0,0)$. So $\langle {\widehat{\alpha }}_j, {\widehat{\alpha }}_i \rangle =0$. \hfill $\square $ \medskip

\noindent \textbf{Corollary 1.2.3b} For every $i,j \in J $ the reflection ${\sigma }_{{\widehat{\alpha }}_i}$ on 
$V$ satisfies  
\begin{displaymath}
{\sigma }_{{\widehat{\alpha }}_i}({\widehat{\alpha }}_j) = 
\left\{ \begin{array}{rl}
-{\widehat{\alpha }}_i, & \mbox{if $j =i$} \\
{\widehat{\alpha }}_j, & \mbox{if $j \ne i$.} \end{array} \right. 
\end{displaymath} 

\noindent \textbf{Proof.} This is immediate from lemma 1.2.3a and the definition of the reflection 
${\sigma }_{{\widehat{\alpha}}_i}$ (\ref{eq-onestar}). \hfill $\square $ \medskip 

\noindent \textbf{Lemma 1.2.3c} For every $i,j,k \in J$ with $i \ne j$ we have 
\begin{equation}
{\sigma }_{{\widehat{\alpha }}_j} \comp {\sigma }_{{\widehat{\alpha }}_i }
({\widehat{\alpha }}_k) = \left\{ \begin{array}{rl}
- {\widehat{\alpha }}_i, & \mbox{if $k=i \, \, \& \, \, j \ne i$} \\
- {\widehat{\alpha }}_j, & \mbox{if $k=j \, \, \& \, \, j \ne i$} \\ 
{\widehat{\alpha }}_k, & \mbox{if $k\ne i \, \, \& \, \, k \ne j$.} \end{array} \right. 
\label{eq-sec1subsec2one}
\end{equation}

\noindent \textbf{Proof.} We compute. If $k \ne i \, \, \& \, \, k \ne j $, then 
\begin{align}
{\sigma }_{{\widehat{\alpha }}_j} \comp {\sigma }_{{\widehat{\alpha }}_i }
({\widehat{\alpha }}_k)  & = {\sigma }_{{\widehat{\alpha }}_j}({\widehat{\alpha }}_k) , 
\quad \mbox{using corollary 1.2.3b} \notag \\
& = {\widehat{\alpha }}_k \quad \mbox{using corollary 1.2.3b.} \notag 
\end{align}
If $k = i \, \, \& \, \, j \ne i$, then 
\begin{align}
{\sigma }_{{\widehat{\alpha }}_j} \comp {\sigma }_{{\widehat{\alpha }}_i }
({\widehat{\alpha }}_i)  & = 
{\sigma }_{{\widehat{\alpha }}_j}(- {\widehat{\alpha }}_i) , 
\quad \mbox{using corollary 1.2.3b} \notag \\
& = -  {\sigma }_{{\widehat{\alpha }}_j}({\widehat{\alpha }}_i)
\quad \mbox{since ${\sigma }_{{\widehat{\alpha }}_j}$ is $\R $-linear } \notag \\
& = - {\widehat{\alpha }}_i \quad \mbox{using $j \ne i$ and corollary 1.2.3b.} \notag 
\end{align}
If $k = j \, \, \& \, \, j \ne i $, then 
\begin{align}
{\sigma }_{{\widehat{\alpha }}_j} \comp {\sigma }_{{\widehat{\alpha }}_i} 
({\widehat{\alpha }}_j)  & = 
{\sigma }_{{\widehat{\alpha }}_j}({\widehat{\alpha }}_j) , 
\quad \mbox{using corollary 1.2.3b.} \notag \\
& = - {\widehat{\alpha }}_j, \quad \mbox{using corollary 1.2.3b.} 
\tag*{$\square $} 
\end{align}

\noindent \textbf{Corollary 1.2.3d} If $i,j \in J $ with $i \ne j $, 
then ${\sigma }_{{\widehat{\alpha }}_i} \comp {\sigma }_{{\widehat{\alpha }}_j} = 
{\sigma }_{{\widehat{\alpha }}_j} \comp {\sigma }_{{\widehat{\alpha }}_i}$. \medskip 

\noindent \textbf{Proof.} For every $k \in J$ we have 
\begin{align}
{\sigma }_{{\widehat{\alpha }}_i} \comp {\sigma }_{{\widehat{\alpha }}_j}
({\widehat{\alpha }}_k) & = \left\{ \begin{array}{rl}
-{\widehat{\alpha }}_i, & \mbox{if $k=i$} \\
- {\widehat{\alpha }}_j, & \mbox{if $k=j$} \\
{\widehat{\alpha }}_k, & \mbox{if $k\ne i \, \, \& \, \, k \ne j$} \end{array} \right.  \notag \\
& = \left\{ \begin{array}{rl}
- {\widehat{\alpha }}_j, & \mbox{if $k=j$} \\
- {\widehat{\alpha }}_i, & \mbox{if $k=i$} \\
{\widehat{\alpha }}_k, & \mbox{if $k\ne i \, \, \& \, \, k \ne j$} \end{array} \right. 
= {\sigma }_{{\widehat{\alpha }}_j} \comp {\sigma }_{{\widehat{\alpha }}_i}
({\widehat{\alpha }}_k). \tag*{$\square $} 
\end{align} 

\noindent \textbf{Corollary 1.2.3e} We have 
\begin{displaymath}
{\sigma }_{{\widehat{\alpha }}_1} \comp {\sigma }_{{\widehat{\alpha }}_2} 
\comp \cdots \comp {\sigma }_{{\widehat{\alpha }}_M} = -{\sigma}_e, 
\end{displaymath}
where ${\sigma}_e$ is the identity mapping of $V$ into itself. \medskip 

\noindent \textbf{Proof.} This follows immediately from (\ref{eq-sec1subsec2one}). 
\hfill $\square $ 

\subsection{A geometric model of $W_{\widehat{\mathcal{R}}}$} 
In this subsection we give a geometric model of the Weyl group $W_{\widehat{\mathcal{R}}}$ of a system 
of roots $\widehat{\mathcal{R}}$ for the nilradical $\widetilde{\mathfrak{n}}$ of a fundamental sandwich algebra. \medskip 

Let ${\widehat{\alpha }}_i \in \widehat{\Pi }$, $i \in J = \{ 1, \ldots , M \} $. For each $i \in J$ there is a unique 
${\widehat{\alpha }}_{j(i)} \in \widehat{\mathcal{R}}$ 
such that ${\widehat{\alpha }}_i + {\widehat{\alpha }}_{j(i)} = \zeta $, that is, 
${\widehat{\alpha }}_{j(i)} = - {\widehat{\alpha }}_i$. For each $i \in J$ let $x_i$ be the root 
vector $X_{{\widehat{\alpha }}_i} \in \widetilde{\mathfrak{n}}$ and $y_i$ be the root vector 
$X_{-{\widehat{\alpha }}_i} 
\in \widetilde{\mathfrak{n}}$. Then 
\begin{displaymath}
[ X_{{\widehat{\alpha }}_i} , X_{{\widehat{\alpha }}_k} ] = \left\{ \begin{array}{rl} 
0, & \mbox{if $k \ne j(i)$} \\
X_{\zeta }, & \mbox{if $k = j(i)$.} \end{array} \right. 
\end{displaymath}
Thus ${\Omega }_{\zeta }(x_i,y_k) = {\nu }_{\zeta }([x_i, y_k]) = {\delta }_{ik}$ is 
a complex valued symplectic form on $Y =L^{+} \oplus L^{-}$, where 
$L^{+} = {\spann }_{\C } \{ x_i \in \widetilde{\mathfrak{n}} \setrule i \in J \} $ and 
$L^{-} = {\spann }_{\C } \{ y_i \in \widetilde{\mathfrak{n}} \setrule i \in J \} $. 
Here ${\nu }_{\zeta }$ is a complex valued linear function on $Z$, which is $1$ on $X_{\zeta }$. $L^{+}$ and $L^{-}$ are Lagrangian subspaces of the symplectic vector space $(Y, 
{\Omega }_{\zeta })$. Thus $Y$ is the ${\Omega }_{\zeta }$ perpendicular direct sum of the 
${\Omega }_{\zeta }$ symplectic planes 
\begin{displaymath}
{\pi }_i = {\spann }_{\C } \{ X_{{\widehat{\alpha }}_i} , X_{-{\widehat{\alpha }}_i} \} 
= {\spann }_{\C } \{ x_i, y_i \} .  
\end{displaymath}
Recall that $\widetilde{\mathfrak{n}} = Y \oplus {\spann }_{\C } \{ X_{\zeta } \} $. \medskip  

For each ${\widehat{\alpha }}_i \in \widehat{\Pi }$, $i \in J$, the reflection 
${\sigma }_{{\widehat{\alpha }}_i}: V \rightarrow V$ gives rise to the $\C $-linear mapping 
\begin{displaymath}
s_{{\widehat{\alpha }}_i}: L^{+} \rightarrow L^{+}: x_k \mapsto \left\{ \begin{array}{rl}
-x_i, & \mbox{if $k =i$} \\
x_k, & \mbox{if $k\ne i$,} \end{array} \right. 
\end{displaymath}
which lifts to the $\C$-linear mapping 
\begin{equation}
\begin{array}{l}
S_i = S_{{\widehat{\alpha }}_i}: Y = L^{+} \oplus L^{-} \rightarrow Y = L^{+} \oplus L^{-}: \\ 
\hspace{.75in} (x_k,y_{\ell }) \mapsto \left\{ \begin{array}{rl} (-x_i,-y_i), & \mbox{if $k=i \, \, \& \, \, \ell =i$} \\
(x_k,y_{\ell }), & \mbox{if $k \ne i $ or $\ell \ne i$.} \end{array} \right. 
\end{array}
\label{eq-nwthreestar}
\end{equation}
So for every $i \in J$ we have $S_i | {\pi }_k = 
${\small $ \left\{ \! \! \begin{array}{rl} -{\mathrm{id}}_{{\pi }_i}, & \mbox{if $k=i$} \\
{\mathrm{id}}_{{\pi }_k}, & \mbox{if $k \ne i $.} \end{array} \right. $} \medskip 

\noindent \textbf{Lemma 1.3.1} For each $i \in J$ the $\C$-linear mapping $S_i$ is a symplectic 
mapping of $(Y, {\Omega }_{\zeta })$ into itself. The following relations hold: 
\begin{displaymath}
\begin{array}{lrl}
1.  & S^2_i & = S_e, \, \mbox{for all $i \in J$}  \\
2. & S_i \comp S_j & = S_j \comp S_i, \, \mbox{for all $i,j \in J$} \\ 
3. & S_1 \comp S_2 \comp \cdots \comp S_M & = - S_e. 
\end{array}
\end{displaymath}
Here $S_e = {\mathrm{id}}_{Y}$. \medskip

\noindent \textbf{Proof.} The proof of the lemma is a straightforward consequence of the definition of 
the mappings $S_i$, $i \in J$. \hfill $\square $ \medskip 
 
\noindent \textbf{Corollary 1.3.1a} Let $\mathcal{W}$ be the group generated by $S_{{\widehat{\alpha }}_i}$, $i \in J$. 
Then the mapping 
\begin{displaymath}
\mu :W_{\widehat{\mathcal{R}}} \rightarrow \mathcal{W}: {\sigma }_{{\widehat{\alpha }}_i} \mapsto S_{{\widehat{\alpha}}_i} 
\end{displaymath} 
is an isomorphism.  

\section{The Weyl group $W_{\widetilde{\mathcal{R}}}$}

Let $\widetilde{\mathfrak{g}} = \mathfrak{g} \oplus \widetilde{\mathfrak{n}}$ be a fundamental 
sandwich algebra. In this section we define the Weyl group $W_{\widetilde{\mathcal{R}}}$ of 
$\widetilde{\mathfrak{g}}$ and show that it is the semidirect product of the Weyl group 
$W_{\mathcal{R}}$ of the simple Lie algebra $\mathfrak{g}$ and the Weyl group $W_{\widehat{\mathcal{R}}}$ 
of the system of roots $\widehat{\mathcal{R}}$ associated to the nilradical $\widetilde{\mathfrak{n}}$ 
of $\widetilde{\mathfrak{g}}$. \medskip 

The system of roots $\widetilde{\mathcal{R}}$ for the fundamental sandwich algebra 
$\widetilde{\mathfrak{g}}$ is the direct sum of two subsystems of 
roots, namely, $\mathcal{R}$, which is a root system of the simple Lie algebra $\mathfrak{g}$ associated 
to the Cartan subalgebra $\mathfrak{h}$, and 
$\widehat{\mathcal{R}}$, which is the system of roots associated to $\mathfrak{h}$ for the nilradical 
$\widetilde{\mathfrak{n}}$ of $\widetilde{\mathfrak{g}}$. Let $\widetilde{U} = U \oplus V $, where 
$U = {\spann }_{\R } \{ \alpha \setrule 
\alpha \in \mathcal{R} \} $ and $V = {\spann }_{\R } \{ \widehat{\alpha } \setrule \widehat{\alpha } \in 
\widehat{\mathcal{R}} \} $. The Weyl group $W_{\mathcal{R}}$ is generated by the reflections 
${\sigma }_{\alpha }$, $\alpha $ a positive root in $\Pi \subseteq \mathcal{R}$, on the 
vector space $U$; whereas the Weyl group $W_{\widehat{\mathcal{R}}}$ is generated by the reflections 
${\sigma }_{\widehat{\alpha }}$, $\widehat{\alpha } \in \widehat{\Pi } \subseteq \widehat{\mathcal{R}}$, 
on the vector space $V$. Let $\widetilde{\Pi } = \Pi \, \amalg \widehat{\Pi }$
be the set of positive roots of the system of roots 
$\widetilde{\mathcal{R}}$. Let  $W_{\widetilde{\mathcal{R}}}$ be the finite group generated by 
linear maps ${\sigma }_{\widetilde{\alpha }}$ of $\widetilde{U}$ into itself such that 
${\sigma }_{\widetilde{\alpha }}|U = {\sigma }_{\alpha }$; while 
${\sigma }_{\widetilde{\alpha }}|V = {\sigma }_{\widehat{\alpha }}$. Here $\widetilde{\alpha } = (\alpha , \widehat{\alpha }) \in \widetilde{\Pi }$. Each ${\sigma }_{\widetilde{\alpha }}$ is an involution on $\widetilde{U}$. The Weyl group $W_{\widehat{\mathcal{R}}}$ 
is a subgroup of $W_{\widetilde{\mathcal{R}}}$ being the image under the injective homomorphism 
\begin{displaymath}
\widetilde{\lambda }: W_{\widehat{\mathcal{R}}} \rightarrow W_{\widetilde{\mathcal{R}}}: 
{\sigma }_{\widehat{\alpha }} \mapsto (1_{W_{\mathcal{R}}}, {\sigma }_{\widehat{\alpha }}). 
\end{displaymath}
The Weyl group $W_{\mathcal{R}}$ is the image of the surjective homomorphism 
\begin{displaymath}
\widetilde{\pi }: W_{\widetilde{\mathcal{R}}} \rightarrow W_{\mathcal{R}}: 
\big( {\sigma }_{\alpha }, {\sigma }_{\widehat{\alpha }} \big) \mapsto {\sigma }_{\alpha }. 
\end{displaymath}
Thus we obtain the sequence
\begin{equation}
1_{W_{\widehat{\mathcal{R}}}} \rightarrow W_{\widehat{\mathcal{R}}} 
\stackrel{\widetilde{\lambda}}{\longrightarrow}  W_{\widetilde{\mathcal{R}}} 
\stackrel{\widetilde{\pi }}{\longrightarrow} W_{\mathcal{R}} \rightarrow 1_{W_{\mathcal{R}}}. 
\label{eq-stardagger}
\end{equation}
By construction of the maps $\widetilde{\lambda }$ and $\widetilde{\pi }$ we have 
$\mathrm{im}\, \widetilde{\lambda } = \ker \widetilde{\pi }$ as sets. Because $\ker \widetilde{\pi }$ is a normal 
subgroup of $W_{\widetilde{\mathcal{R}}}$, the group $\widetilde{\lambda }(W_{\widehat{\mathcal{R}}})$ 
must be a normal subgroup of $W_{\widetilde{\mathcal{R}}}$. It is, because $W_{\widehat{\mathcal{R}}}$ is an 
abelian group and hence $\widetilde{\lambda }(W_{\widehat{\mathcal{R}}})$ is also. Thus the sequence 
(\ref{eq-stardagger}) is an exact sequence of groups. \medskip 

We now prove \medskip 

\noindent \textbf{Theorem 2.1} The Weyl group $W_{\widetilde{\mathcal{R}}}$ of the system of roots 
$\widetilde{\mathcal{R}}$ associated to the fundamental sandwich algebra $\widetilde{\mathfrak{g}} = 
\mathfrak{g} \oplus \widetilde{\mathfrak{n}}$ is the semidirect product $W_{\widehat{\mathcal{R}}} {\rtimes}_{\varphi } W_{\mathcal{R}}$ of the abelian Weyl group $W_{\widehat{\mathcal{R}}} = \mathcal{W}$ 
of the system of roots $\widehat{\mathcal{R}}$ associated to the nilpotent radical $\widetilde{\mathfrak{n}}$ of 
$\widetilde{\mathfrak{g}}$ and the Weyl group $W_{\mathcal{R}}$ associated to the root system $\mathcal{R}$ of the simple Lie algebra $\mathfrak{g}$. Here $\varphi : W_{\mathcal{R}} \rightarrow \mathrm{Aut}(\mathcal{W})$ is a group homomorphism from $W_{\mathcal{R}}$ into the group $\mathrm{Aut}(\mathcal{W})$ of automorphisms of $\mathcal{W}$, which is involved in defining the group multiplication in $W_{\widehat{\mathcal{R}}} {\rtimes}_{\varphi } W_{\mathcal{R}}$, see equation (\ref{eq-ddagger}) below. \medskip 

\noindent \textbf{Proof.} To prove the theorem we need to construct the homomorphism $\varphi $, 
which is a consequence of the following lemmas. \medskip 

We start with the following construction. Let $W_{\mathcal{R}}$ be the Weyl group 
associated to root system $\mathcal{R}$ of the simple Lie algebra 
$\mathfrak{g}$ with Cartan subalgebra $\mathfrak{h}$. Suppose that $\mathfrak{g}$ is a 
subalgebra of a simple Lie algebra $\underline{\mathfrak{g}}$ with root system $\underline{\mathcal{R}}$ 
associated to a Cartan subalgebra $\underline{\mathfrak{h}}$. Suppose that the 
Cartan subalgebra $\mathfrak{h}$ is \emph{aligned} with the Cartan subalgebra $\underline{\mathfrak{h}}$, 
that is, there is a vector $\underline{\widetilde{H}} \in \underline{\mathfrak{h}}$ such that 
$\mathcal{R} = \{ \underline{\alpha } \in \underline{\mathcal{R}} \setrule 
\underline{\alpha }(\underline{\widetilde{H}}) =0 \} $. Let ${\mathcal{R}}^{-} 
=\{ \underline{\alpha } \in \underline{\mathcal{R}} \setrule 
\underline{\alpha }(\underline{\widetilde{H}}) <0 \} $. 
From \cite{cushman17b} recall that the system of roots $\widehat{\mathcal{R}}$ of the nilradical $\widetilde{\mathfrak{n}}$ of $\widetilde{\mathfrak{g}}$ is the collection of linear functions $\widehat{\alpha }$ on $\mathfrak{h}$ such that there is a linear function $\underline{\alpha }$ in ${\mathcal{R}}^{-}$ whose restriction to $\mathfrak{h}$ is $\widehat{\alpha }$.  \medskip 
 
\noindent \textbf{Lemma 2.1a} The Weyl group $W_{\mathcal{R}}$ acts on ${\mathcal{R}}^{-}$. \medskip 

\noindent \textbf{Proof}. Suppose that $\alpha \in \mathcal{R}$. Then there is an $\underline{\alpha } \in 
\underline{\mathcal{R}}$ such that $\alpha = \underline{\alpha }|\mathfrak{h}$. Let 
${\sigma }_{\underline{\alpha }} \in W_{\underline{\mathcal{R}} }$ be the reflection in 
$\underline{V} = {\spann }_{\R}\{ \underline{\beta } \in \underline{\mathcal{R}} \}$ corresponding to the 
root $\underline{\alpha }$. Then ${\sigma }_{\alpha } = {\sigma }_{\underline{\alpha }}|V$, where 
$V = {\spann }_{\R }\{ \widehat{\alpha } \in \widehat{\mathcal{R}} \} $. For each 
$\underline{\beta } \in {\mathcal{R}}^{-} \subseteq \underline{\mathcal{R}}$ we have 
${\sigma }_{\underline{\alpha}}(\underline{\beta }) \in \underline{\mathcal{R}} \subseteq {\underline{\mathfrak{h}}}^{\ast}$, since $\underline{\mathcal{R}}$ is a root system for the simple 
Lie algebra $\underline{\mathfrak{g}}$. So
\begin{align}
{\sigma }_{\underline{\alpha }}(\underline{\beta })(\underline{\widetilde{H}}) & = 
\underline{\beta }(\underline{\widetilde{H}}) 
- \langle \underline{\beta }, \underline{\alpha } \rangle \underline{\alpha }(\underline{\widetilde{H}}) 
= \underline{\beta }(\underline{\widetilde{H}}), \, \, \, \mbox{since $\alpha \in \mathcal{R}$} \notag \\
& < 0, \, \, \, \mbox{since $\underline{\beta } \in {\mathcal{R}}^{-}$.} \notag 
\end{align}  
Therefore every reflection ${\sigma }_{\alpha } \in W_{\mathcal{R}}$ induces a linear mapping 
\begin{equation}
s_{\alpha}: \widehat{\mathcal{R}} \subseteq V \rightarrow \widehat{\mathcal{R}} \subseteq V:
\widehat{\beta } = \underline{\beta } |\mathfrak{h} \mapsto s_{\alpha }(\widehat{\beta }) = 
{\sigma }_{\underline{\alpha }}(\underline{\beta })|\mathfrak{h},  
\label{eq-nwstar}
\end{equation}
where $\underline{\beta } \in {\mathcal{R}}^{-}$. The mapping $s_{\alpha }$ (\ref{eq-nwstar}) is well defined. To see this suppose that $\widehat{\beta } = \underline{\gamma}|\mathfrak{h}$ for some $\underline{\gamma } \in 
{\mathcal{R}}^{-}$. Since the reflection ${\sigma }_{\underline{\alpha }}$ induces a complex linear mapping 
${\sigma }_{\underline{\alpha }}:{\underline{\mathfrak{h}}}^{\ast } \rightarrow {\underline{\mathfrak{h}}}^{\ast }$, 
we get  
\begin{align}
{\sigma }_{\underline{\alpha }}(\underline{\gamma })|\mathfrak{h} & =
\underline{\gamma }({\sigma }^T_{\underline{\alpha } }|\mathfrak{h}) 
= \underline{\beta }({\sigma }^T_{\underline{\alpha }}|\mathfrak{h}), \, \, \, 
\mbox{because $\underline{\gamma }|\mathfrak{h} = \underline{\beta }|\mathfrak{h}$} \notag \\
& = {\sigma }_{\underline{\alpha}}(\underline{\beta })|\mathfrak{h}, \quad 
\mbox{because ${\sigma }^T_{\underline{\alpha }}|\mathfrak{h}$ maps $\mathfrak{h}$ into itself.} 
\tag*{$\square $}
\end{align}

We begin the construction of the semidirect product by looking at the finer structure of the 
linear mappings $s_{\alpha }$ (\ref{eq-nwstar}) on the vector space $V$. \medskip  

\noindent \textbf{Corollary 2.1b} The mapping $s_{\alpha }$, $\alpha \in \mathcal{R}$, is an involution. \medskip 

\noindent \textbf{Proof.} This follows because 
\begin{align}
s_{\alpha }\big( s_{\alpha }(\widehat{\beta }) \big) & = 
s_{\alpha }\big( {\sigma }_{\underline{\alpha }}(\underline{\beta }|\mathfrak{h}) \big) = 
{\sigma }_{\underline{\alpha }}\big( {\sigma }_{\underline{\alpha }}(\underline{\beta }) \big) |\mathfrak{h} 
= \underline{\beta }|\mathfrak{h} = \widehat{\beta }. \tag*{$\square $} 
\end{align}

Recall that the vector space 
$V = {\spann }_{\R } \{ \widehat{\alpha } \setrule \widehat{\alpha } \in \widehat{\Pi} \} $, 
has $\widehat{\mathcal{R}}$ as its abelian group of vectors. For each positive root $\alpha \in \mathcal{R}$ 
the Weyl group $W_{\mathcal{R}}$ acts linearly on $V$ via the involution 
$s_{\alpha }$ induced from the reflection ${\sigma }_{\alpha }$. For every 
${\widehat{\alpha}}_j \in \widehat{\Pi }$, $i \in J$, we have 
\begin{align*}
s_{\alpha }(\zeta ) & = s_{\alpha } \big({\widehat{\alpha }}_j +(- {\widehat{\alpha }}_j) \big)  
= s_{\alpha }({\widehat{\alpha }}_j) +s_{\alpha }(-{\widehat{\alpha }}_j)   
 = s_{\alpha }({\widehat{\alpha }}_j) - s_{\alpha }({\widehat{\alpha }}_j) = \zeta , 
\end{align*}
since $s_{\alpha }:V \rightarrow V$ is a linear map, 
${\widehat{\alpha }}_j \in \widehat{\mathcal{R}}$ implies that $s_{\alpha }({\widehat{\alpha }}_j)$ lies in 
$\widehat{\mathcal{R}}$ by lemma 2.1a, and $\widehat{\mathcal{R}} + (-\widehat{\mathcal{R}}) = \{ \zeta  \} $. 
Thus $s_{\alpha }$ maps $\widehat{\mathcal{R}} $ bijectively to itself. \medskip 

Since $\widehat{\mathcal{R}} = \widehat{\Pi } \amalg (-\widehat{\Pi })$, for each $i \in J = \{ 1, \ldots , M \} $ the image of the root ${\widehat{\alpha }}_i$ in $\widehat{\mathcal{R}} $ under the mapping $s_{\alpha }$ lies in either $\widehat{\Pi }$ or $-\widehat{\Pi }$ but not both. Define a permutation ${\tau }_{s_{\alpha }}$ on 
for each $i \in J$ by 
\begin{displaymath}
{\widehat{\alpha }}_{{\tau }_{s_{\alpha }}(i)} = 
\left\{ \! \! \! \begin{array}{rl} 
s_{\alpha }({\widehat{\alpha }}_i), & \mbox{if $s_{\alpha }({\widehat{\alpha }}_i) \in \widehat{\Pi }$} \\
-s_{\alpha }({\widehat{\alpha }}_i), & \mbox{if $s_{\alpha }({\widehat{\alpha }}_i) \in -\widehat{\Pi }$.} 
\end{array} \right. 
\end{displaymath} 

\noindent \textbf{Lemma 2.1c} If $\alpha $, $\beta \in \Pi \subseteq \mathcal{R}$, then 
${\tau }_{s_{\alpha } \comp s_{\beta }} = {\tau }_{s_{\alpha }} \comp {\tau }_{s_{\beta }} $. In other 
words, the mapping 
\begin{displaymath} 
W_{\mathcal{R}} \rightarrow {\mathrm{Perm}}_J: {\sigma }_{\alpha } \mapsto {\tau }_{s_{\alpha }}
\end{displaymath}
is a group homomorphism. \medskip 

\noindent \textbf{Proof.} This is an immediate consequence of the definition of the mapping 
$s_{\alpha }$ (\ref{eq-nwstar}) and the permutation ${\tau }_{s_{\alpha }}$. \hfill $\square $ \medskip 

The linear mapping 
\begin{displaymath}
t_{{\tau }_{s_{\alpha }}}: L^{+} \rightarrow L^{+}: X_{{\widehat{\alpha }}_i} \mapsto 
X_{{\widehat{\alpha }}_{{\tau }_{s_{\alpha }}(i)}}, \, \, \, \mbox{for each $i \in J$} 
\end{displaymath}
is well defined because ${\widehat{\alpha }}_{{\tau }_{s_{\alpha }}(i)} \in \widehat{\Pi }$ for every 
$i \in J $. Lift the mapping $t_{{\tau }_{s_{\alpha }}}$ to a linear symplectic mapping 
of $(Y, {\Omega }_{\zeta })$ into itself given by 
\begin{displaymath}
T_{{\tau }_{s_{\alpha }}}: Y = L^{+} \oplus L^{-} \rightarrow Y= L^{+} \oplus L^{-}: 
\big( X_{{\widehat{\alpha }}_i}, X_{-{\widehat{\alpha }}_j} \big) 
\longmapsto \big( X_{{\widehat{\alpha }}_{{\tau }_{s_{\alpha }}(i)}}, X_{-{\widehat{\alpha }}_{{\tau }_{s_{\alpha }}(j)}} \big) .  
\end{displaymath}
From its construction we see that the linear mapping $T_{{\tau }_{s_{\alpha }}}$ sends the 
${\Omega }_{\zeta }$ symplectic $2$-plane ${\pi }_i = {\spann }_{\R }\{ X_{{\widehat{\alpha }}_i}, 
X_{-{\widehat{\alpha }}_i} \} $ onto the ${\Omega}_{\zeta }$ symplectic $2$-plane 
${\pi }_{{\tau }_{s_{\alpha }}(i)} = {\spann }_{\R }\{ X_{{\widehat{\alpha }}_{{\tau }_{s_{\alpha }}(i)}}, 
X_{-{\widehat{\alpha }}_{{\tau }_{s_{\alpha }}(i)}} \} $ for every $i \in J$. \medskip 

\noindent \textbf{Lemma 2.1d} For every ${\sigma }_{\alpha }$, ${\sigma }_{\beta } \in W_{\mathcal{R}}$ 
we have $T_{{\tau }_{s_{\alpha } \comp s_{\beta } }} = T_{{\tau }_{s_{\alpha }}} \comp T_{{\tau }_{s_{\beta }}}$. 
\medskip 

\noindent \textbf{Proof.} This follows from the construction of the mapping $T_{{\tau }_{s_{\alpha }}}$. 
\hfill $\square $ \medskip 

\noindent \textbf{Lemma 2.1e} For every $i \in J $ we have 
\begin{align}
T_{{\tau }_{s_{\alpha }}} \comp S_{{\widehat{\alpha }}_i} \comp T^{-1}_{{\tau }_{s_{\alpha }}} & = 
S_{{\widehat{\alpha }}_{{\tau }_{s_{\alpha }}(i)}}. 
\label{eq-star}
\end{align} 

\noindent \textbf{Proof.} Consider the $2$-plane ${\pi }_{{\tau }_{s_{\alpha }}(i)}$ then 
$T^{-1}_{{\tau }_{s_{\alpha }}}\big( {\pi }_{{\tau }_{s_{\alpha }}(i)} \big) = {\pi }_i$; while 
$S_{{\widehat{\alpha }}_i}|{\pi }_i = -{\mathrm{id}}_{{\pi }_i}$. So $T_{{\tau }_{s_{\alpha }}} \comp S_{{\widehat{\alpha }}_i} \comp T^{-1}_{{\tau }_{s_{\alpha }}} ( {\pi }_{{\tau }_{s_{\alpha }}(i)} ) = 
- {\pi }_{{\tau }_{s_{\alpha }}(i)}$. If $j \ne i$ then $T^{-1}_{{\tau }_{s_{\alpha }}}\big( {\pi }_{{\tau }_{s_{\alpha }}(j)} \big) = 
{\pi }_j$; while $S_{{\widehat{\alpha }}_j}|{\pi }_j = {\mathrm{id}}_{{\pi }_j}$. So 
$T_{{\tau }_{s_{\alpha }}} \comp S_{{\widehat{\alpha }}_i} \comp T^{-1}_{{\tau }_{s_{\alpha }}} 
( {\pi }_{{\tau }_{s_{\alpha }}(j)} ) = {\pi }_{{\tau }_{s_{\alpha }}(j)} $. Consequently, (\ref{eq-star}) holds 
by definition of $S_{{\widehat{\alpha}}_i}$ (\ref{eq-nwthreestar}).  \hfill $\square $\medskip 

\noindent \textbf{Lemma 2.1f} The mapping 
\begin{equation}
{\varphi }_{s_{\alpha }}: \mathcal{W} \rightarrow \mathcal{W}: 
S_{{\widehat{\alpha }}_i} \mapsto T_{{\tau }_{s_{\alpha }}} \comp S_{{\widehat{\alpha }}_i} \comp 
T^{-1}_{s_{\alpha }}
\label{eq-ddaggernw}
\end{equation}
is an isomorphism. \medskip 

\noindent \textbf{Proof.} Since $\{ S_{{\widehat{\alpha }}_i}, \, i \in J \} $ generates 
$\mathcal{W} = W_{\widehat{\mathcal{R}}}$, every element $w \in \mathcal{W}$ may be written as $w = S_{{\widehat{\alpha }}_{i_1}} \comp \cdots \comp S_{{\widehat{\alpha }}_{i_{\ell }}}$. Thus 
${\varphi }_{s_{\alpha }}$ maps $\mathcal{W}$ into 
itself, because 
\begin{align}
T_{{\tau }_{s_{\alpha }}}  w \, T^{-1}_{s_{\alpha }} & = 
T_{{\tau }_{s_{\alpha }}} \comp \big( S_{{\widehat{\alpha }}_{i_1}} \comp 
\cdots \comp S_{{\widehat{\alpha }}_{i_{\ell }}} \big) \comp  T^{-1}_{s_{\alpha }}  \notag \\
& = \big( T_{{\tau }_{s_{\alpha }}} \comp S_{{\widehat{\alpha }}_{i_1}} \comp T^{-1}_{s_{\alpha }} \big) 
\comp \cdots \comp \big( T_{{\tau }_{s_{\alpha }}} \comp S_{{\widehat{\alpha }}_{i_{\ell }}} \comp 
T^{-1}_{s_{\alpha }} \big) 
\notag \\
& = S_{{\widehat{\alpha }}_{{\tau }_{s_{\alpha }}(i_1)}} \comp 
\cdots \comp S_{{\widehat{\alpha }}_{{\tau }_{s_{\alpha }}(i_{\ell }) }} \in \mathcal{W}. \notag 
\end{align} 
Clearly for $w$, $w' \in \mathcal{W}$ we have 
$T_{{\tau }_{s_{\alpha }}}  (ww')\, T^{-1}_{s_{\alpha }} = 
\big( T_{{\tau }_{s_{\alpha }}}  w \, T^{-1}_{s_{\alpha }}\big) 
\big( T_{{\tau }_{s_{\alpha }}}  w' \, T^{-1}_{s_{\alpha }} \big) $, 
that is, the map ${\varphi }_{s_{\alpha }}$ (\ref{eq-ddaggernw}) is a homomorphism. Since 
${\varphi }^{-1}_{s_{\alpha }} = 
{\varphi }_{s^{-1}_{\alpha }}$, the map ${\varphi }_{s_{\alpha }}$ is a isomorphism, that is, 
an automorphism of $\mathcal{W}$. \hfill $\square $ \medskip  

\noindent \textbf{Lemma 2.1g} The map 
\begin{equation}
\varphi : W_{\mathcal{R}} \rightarrow \mathrm{Aut}(\mathcal{W}): {\sigma }_{\alpha } \mapsto 
{\varphi }_{s_{\alpha }}
\label{eq-tdagger}
\end{equation}
is a group homomorphism. \medskip 

\noindent \textbf{Proof.} We need only show that for every ${\sigma }_{\alpha }$, ${\sigma }_{\beta } \in 
\mathcal{W}$ we have ${\varphi }_{{\sigma }_{\alpha } \comp {\sigma }_{\beta }} = 
{\varphi }_{{\sigma }_{\alpha }} \comp {\varphi }_{{\sigma }_{\beta }} $. For every $w \in \mathcal{W}$ we get 
\begin{align}
{\varphi }_{{\sigma }_{\alpha }} \big( {\varphi }_{{\sigma }_{\beta}}(w) \big) & = 
{\varphi }_{{\sigma }_{\alpha }}\big( T_{{\tau }_{s_{\beta }}} w T^{-1}_{{\tau }_{s_{\beta }}} \big) 
 = T_{{\tau }_{s_{\alpha }}}\big( T_{{\tau }_{s_{\beta }}} w \, T^{-1}_{{\tau }_{s_{\beta }}} \big) 
T^{-1}_{{\tau }_{s_{\alpha }}} \notag \\
&\hspace{-.75in} = (T_{{\tau }_{s_{\alpha }}}T_{{\tau }_{s_{\beta }}}) w\, (T_{{\tau }_{s_{\alpha }}}T_{{\tau }_{s_{\beta }}})^{-1}
= T_{{\tau }_{s_{\alpha }} \comp {\tau }_{s_{\beta }}} w \, T^{-1}_{{\tau }_{s_{\alpha }} \comp {\tau }_{s_{\beta }}} 
= {\varphi }_{{\sigma }_{\alpha } \comp {\sigma }_{\beta }}(w). \tag*{$\square $}
\end{align}

We now use the map $\varphi $ (\ref{eq-tdagger}) to define a multiplication {\tiny \raisebox{1pt}{$\bullet $}} on 
$\mathcal{W} \times W_{\mathcal{R}}$, which gives the semidirect product 
$\mathcal{W} {\rtimes}_{\varphi } W_{\mathcal{R}}$. Define the multiplication {\tiny \raisebox{1pt}{$\bullet $}} 
by 
\begin{equation}
(w, {\sigma }_{\alpha }) \mbox{\tiny \raisebox{1pt}{$\bullet $}} (w', {\sigma }_{\beta }) = 
(w {\varphi }_{{\sigma }_{\alpha }}(w'), {\sigma }_{\alpha } \comp {\sigma }_{\beta } ) . 
\label{eq-ddagger}
\end{equation}

Next we show that the multiplication operation 
{\tiny \raisebox{1pt}{$\bullet $}} (\ref{eq-ddagger}) turns $\mathcal{W} \times W_{\mathcal{R}}$ into a 
group $\mathcal{W} {\rtimes}_{\varphi } W_{\mathcal{R}}$. The following calculation shows that $(1_{\mathcal{W}}, 1_{W_{\mathcal{R}}})$ is 
the identity element of $\mathcal{W} {\rtimes}_{\varphi } W_{\mathcal{R}}$ because  
\begin{align} 
(1_{\mathcal{W}}, 1_{W_{\mathcal{R}}} ) \mbox{\tiny \raisebox{1pt}{$\bullet $}} (w,{\sigma }_{\alpha }) & = 
\big( 1_{\mathcal{W}} {\varphi }_{1_{W_{\mathcal{R}}}}(w), 1_{W_{\mathcal{R}}} {\sigma }_{\alpha} \big) 
= ({\varphi }_{1_{\mathcal{W}}} (w), 1_{W_{\mathcal{R}}} {\sigma }_{\alpha } ) 
\notag \\
& = \big(w{\varphi }_{{\sigma }_{\alpha }}(1_{\mathcal{W}}),{\sigma }_{\alpha }1_{W_{\mathcal{R}}} \big) ) = 
(w,{\sigma }_{\alpha }) \mbox{\tiny \raisebox{1pt}{$\bullet $}} (1_{\mathcal{W}}, 1_{W_{\mathcal{R}}} ). \notag 
\end{align}
The inverse of $(w,{\sigma }_{\alpha })$ is 
$\big( {\varphi }_{({\sigma }_{\alpha })^{-1}}(w^{-1}), ({\sigma }_{\alpha })^{-1} \big) $, since 
\begin{align}
(w,{\sigma }_{\alpha } ) \mbox{\tiny \raisebox{1pt}{$\bullet $}} \big( {\varphi }_{({\sigma }_{\alpha })^{-1}}(w^{-1}), 
({\sigma }_{\alpha })^{-1} \big) & = 
\big( w {\varphi }_{{\sigma }_{\alpha }}({\varphi }_{({\sigma }_{\alpha })^{-1}}(w^{-1})), 
{\sigma }_{\alpha} ({\sigma }_{\alpha })^{-1} \big) \notag \\
& = \big( w w^{-1}, {\sigma }_{\alpha } ({\sigma }_{\alpha })^{-1} \big) = (1_{\mathcal{W}},1_{W_{\mathcal{R}}} ) \notag 
\end{align}
and 
\begin{align}
\big( {\varphi }_{({\sigma }_{\alpha})^{-1}}(w^{-1}), ({\sigma }_{\alpha })^{-1} \big) \mbox{\tiny \raisebox{1pt}{$\bullet $}} 
(w,{\sigma }_{\alpha }) & = 
\big( {\varphi }_{({\sigma }_{\alpha })^{-1}}(w^{-1}) {\varphi }_{({\sigma }_{\alpha })^{-1}}(w), ({\sigma }_{\alpha })^{-1} 
{\sigma }_{\alpha } \big) 
\notag \\ 
& = \big( {\varphi}_{({\sigma }_{\alpha })^{-1}}(w^{-1}w), ({\sigma }_{\alpha })^{-1}{\sigma }_{\alpha } \big) = 
(1_{\mathcal{W}}, 1_{W_{\mathcal{R}}}). \notag 
\end{align}
The multiplication $\mbox{\tiny \raisebox{1pt}{$\bullet $}}$ is associative, because  
\begin{align}
\big( (w,{{\sigma}_{\alpha}})  \mbox{\tiny \raisebox{1pt}{$\bullet $}} (w',{\sigma}'_{\alpha}) \big) \mbox{\tiny \raisebox{1pt}{$\bullet $}} (w'', {\sigma}''_{\alpha}) & = 
\big( w {\varphi }_{{\sigma}_{\alpha}}(w'), {\sigma}_{\alpha} {\sigma}'_{\alpha} \big) 
\mbox{\tiny \raisebox{1pt}{$\bullet $}} (w'', {\sigma}''_{\alpha} ) 
\notag \\
& \hspace{-1.75in} = \big( w {\varphi }_{{\sigma}_{\alpha}}(w') {\varphi }_{{\sigma}_{\alpha} {\sigma}'_{\alpha}}(w''), 
{\sigma}_{\alpha} {\sigma}'_{\alpha} {\sigma}''_{\alpha} \big) 
= \big( w {\varphi }_{{\sigma}_{\alpha}} \big( w' {\varphi }_{{\sigma}'_{\alpha}}(w'') \big), 
{\sigma}_{\alpha} {\sigma}'_{\alpha} {\sigma}''_{\alpha} \big) 
\notag \\
& \hspace{-1.75in} = (w,{\sigma}_{\alpha}) \mbox{\tiny \raisebox{1pt}{$\bullet $}} 
\big( w' {\varphi }_{{\sigma}'_{\alpha}}(w''), {\sigma}'_{\alpha} {\sigma}''_{\alpha} \big) 
= (w, {\sigma}_{\alpha} )  \mbox{\tiny \raisebox{1pt}{$\bullet $}} \big( (w', {\sigma}'_{\alpha} ) 
\mbox{\tiny \raisebox{1pt}{$\bullet $}} (w'', {\sigma}''_{\alpha} ) \big) . \notag 
\end{align}
Thus the pair $(\mathcal{W} {\rtimes}_{\varphi }  W_{\mathcal{R}} , \mbox{\tiny \raisebox{1pt}{$\bullet $}}  )$ 
is a group. \medskip  

Consider the homomorphisms 
\begin{displaymath}
\widetilde{\lambda }: \mathcal{W} \rightarrow \mathcal{W} {\rtimes}_{\varphi } W_{\mathcal{R}}: 
w \mapsto (w, 1_{W_{\mathcal{R}}})
\end{displaymath}
and 
\begin{displaymath}
\widetilde{\pi} : \mathcal{W} {\rtimes}_{\varphi } W_{\mathcal{R}} \rightarrow W_{\mathcal{R}}: (w, {\sigma }_{\alpha }) 
\mapsto {\sigma }_{\alpha }. 
\end{displaymath}
Together they give the exact sequence
\begin{align}
& 1_{\mathcal{W}} \rightarrow \mathcal{W} \stackrel{\widetilde{\lambda }}{\longrightarrow} \mathcal{W} {\rtimes}_{\varphi } W_{\mathcal{R}} \stackrel{\widetilde{\pi} }{\longrightarrow} W_{\mathcal{R}} \rightarrow 1_{W_{\mathcal{R}}} , 
\label{eq-doublestar}
\end{align}
which is the same as the exact sequence (\ref{eq-stardagger}). Fix $w \in \mathcal{W}$. Define the mapping 
\begin{displaymath}
\gamma : W_{\mathcal{R}} \rightarrow \mathcal{W} {\rtimes}_{\varphi } W_{\mathcal{R}}: 
{\sigma }_{\alpha } \mapsto \widetilde{\lambda} (w^{-1}) \mbox{\tiny \raisebox{1pt}{$\bullet $}} 
( w , {\sigma }_{\alpha } ) .
\end{displaymath}
The following argument shows that $\gamma $ is well defined. We have 
\begin{align}
\gamma ({\sigma }_{\alpha }) & = 
\widetilde{\lambda } (w^{-1}) \mbox{\tiny \raisebox{1pt}{$\bullet $}} ( w , {\sigma }_{\alpha } ) \notag \\ 
& = \big( w^{-1}, 1_{W_{\mathcal{R}}} \big) \mbox{\tiny \raisebox{1pt}{$\bullet $}} ( w, {\sigma }_{\alpha } )  
= ( w^{-1} w , 1_{W_{\mathcal{R}}} \comp {\sigma }_{\alpha } ) = (1_{\mathcal{W}}, {\sigma }_{\alpha }), \notag 
\end{align}
which does not depend on the choice of $w \in \mathcal{W}$. Clearly $\gamma $ is a group homomorphism. 
Moreover, $\widetilde{\pi } \comp \gamma = {\mathrm{id}}_{W_{\mathcal{R}}}$. Thus $\gamma $ splits 
the sequence $(\ref{eq-doublestar})$. This completes the proof that the Weyl group 
$W_{\widetilde{\mathcal{R}}}$ is the semidirect product $W_{\widehat{\mathcal{R}}} {\rtimes}_{\varphi } W_{\mathcal{R}}$, where $W_{\widehat{\mathcal{R}}} = \mathcal{W}$. \hfill $\square $


\begin{thebibliography}{99}

\bibitem{cushman17b} Cushman, R., Systems of roots, \texttt{arXiv.1708.02568}. 

\end{thebibliography}
\end{document}